\documentclass[11pt,bezier]{article}
\usepackage{amsmath,amssymb,amsfonts,euscript,graphicx}

\newtheorem{The}{Theorem}[section]
\newtheorem{Pro}[The]{Proposition}
\newtheorem{Lem}[The]{Lemma}
\newtheorem{Cor}[The]{Corollary}
\newtheorem{Def}[The]{Definition}

\newtheorem{Examp}[The]{Example}

\title{\bf \Large  Modules Satisfying the Prime Radical Condition\\ and a Sheaf Construction
for Modules I\thanks {The research
 of first author   was in part supported by
a grant from IPM (No. 90160034).}
\thanks {2010{ \it Mathematics Subject Classification}: 13C13, 13C99, 13A99, 14A25.
 }}
\author{{\normalsize \bf { M. Behboodi$^{{\rm a,b}}$\thanks{Corresponding
 author.} and M. Sabzevari$^{{\rm a}}$}}\\
 {\small{ $^{\rm a}$Department of Mathematical Sciences, Isfahan University of Technology}}\vspace{-1mm}\\
  {\small{ P.O.Box: 84156-83111, Isfahan, Iran}}\\
 {\small{ $^{\rm b}$School of Mathematics, Institute for Research in Fundamental Sciences
 (IPM)}}\vspace{-1mm}\\ {\small{ P.O.Box: 19395-5746, Tehran, Iran}}\vspace{-1mm}\\
{\footnotesize{$\mathsf{mbehbood@cc.iut.ac.ir}$}}\vspace{-1mm}\\
{\footnotesize{$\mathsf{sabzevari@math.iut.ac.ir}$}}}
\date{}
\begin{document}
\maketitle
\begin{abstract}
{\small \noindent The purpose of this paper and its sequel, is to
introduce a new class of modules over a commutative ring $R$,
called $\mathbb{P}$-radical modules (modules $M$ satisfying the
prime radical condition ``$(\sqrt[p]{{\cal{P}}M}:M)={\cal{P}}$''
for every prime ideal ${\cal{P}}\supseteq {\rm Ann}(M)$, where
$\sqrt[p]{{\cal{P}}M}$ is the intersection of all prime submodules
of $M$ containing ${\cal{P}}M$). This class contains the family of
primeful modules properly. This yields that over any ring all free
modules and all finitely generated modules lie in the class of
$\mathbb{P}$-radical modules. Also, we show that if $R$ is a
domain  (or a Noetherian ring), then   all projective modules are
$\mathbb{P}$-radical. In particular, if $R$ is an Artinian ring,
then all $R$-modules are $\mathbb{P}$-radical and the converse is
also true when $R$ is a Noetherian ring.  Also an $R$-module $M$
is called $\mathbb{M}$-radical if
$(\sqrt[p]{{\cal{M}}M}:M)={\cal{M}}$; for every maximal ideal
${\cal{M}}\supseteq {\rm Ann}(M)$.  We show that the two concepts
$\mathbb{P}$-radical and $\mathbb{M}$-radical are equivalent for
all $R$-modules if and only if $R$ is a Hilbert ring. Semisimple
$\mathbb{P}$-radical ($\mathbb{M}$-radical) modules are also
characterized. In Part II we shall continue the study of this
construction,  and as an application, we show that the sheaf
theory of spectrum of $\mathbb{P}$-radical modules (with the
Zariski topology) resembles to that of rings.}
\end{abstract}

\section{Introduction}
All rings in this article are commutative with identity and
modules are unital. For a ring $R$ we denote by $dim(R)$ the
classical Krull dimension of $R$ and for a  submodule $N$ of an
$R$-module $M$ we denote the annihilator of the factor module
$M/N$  by $(N:M)$, i.e., $(N:M)=\{r\in R~|~rM\subseteq N\}$.  We
call $M$ faithful if $(0:M)=0.$

The theory sheaf of rings (with the Zariski toplogy)  on the
spectrum of prime ideals of a commutative ring is one of the main
tools in Algebraic Geometry. Recall that the {\it spectrum } {\rm
Spec}($R$) of a ring $R$ consists of all prime ideals of $R$ and
is non-empty. For each
 ideal $I$ of $R$, we set ${V(I)=\{{\cal{P}}\in {\rm Spec}(R) :
I\subseteq {\cal{P}}\}}$. Then the sets $V(I)$, where $I$ is an
ideal of $R$,  satisfy the axioms for the closed sets of a
topology on {\rm Spec}($R$), called the {\it Zariski topology} of
$R$. The distinguished open sets of {\rm Spec}($R$) are the open
sets of the form $D(f)=\{{\cal{P}}\in {\rm Spec}(R) : f\not\in
{\cal{P}}\}={\rm Spec}(R)\setminus V(f)$, where $V(f)=V(Rf)$.
These sets form a basis for the Zariski topology on {\rm
Spec}$(R)$ (see for examples, Atiyah and Macdonald \cite{Atiyah}
and Hartshorne \cite{Hartshorne}).

We recall that a {\it sheaf of rings} ${\cal{O}}_X$ on a
topological space $X$ is an assignment of a ring ${\cal{O}}_X(U)$
to each open set $U$ in $X$, together with, for each inclusion
 $U\subseteq V$ a {\it restriction homomorphism} $res_{V,U} : {\cal{O}}_X(V)\rightarrow {\cal{O}}_X(U)$, subject to the
following conditions:\\
\indent (i) ${\cal{O}}_X(\emptyset)=0$.\\
 \indent (ii) $res_{U,U}=id_U$.\\
  \indent(iii)  If $U\subseteq  V\subseteq W$, then $res_{V,U}\circ res_{W,V}=res_{W,U}$.\\
  \indent (iv) For each open cover $\{U_\alpha\}$ of $U\subseteq X$  and for each
collection of elements ${f_\alpha\in{\cal{O}}_X(U_\alpha)}$
\indent\indent such that for all $\alpha,~\beta$ if
$res_{U_\alpha,U_\alpha\cap
U_\beta}(f_\alpha)=res_{U_\beta,U_\alpha\cap U_\beta}(f_\beta)$,
then there is a unique \indent\indent $f\in {\cal{O}}_X(U)$ such
that for all $\alpha$, $f_\alpha=res_{U,U_\alpha}(f)$.

We think of the elements of ${\cal{O}}_X(U)$ as ``functions"
defined on $U$. The restriction homomorphisms correspond to
restricting a function on a big open set to a smaller one.

It is well-known that for any commutative ring $R$, there is a
sheaf of rings on {\rm Spec}$(R)$, called the {\it structure
sheaf}, denoted by ${\cal{O}}_{{\rm Spec}(R)}$, defined as
follows: for each prime ideal ${\cal{P}}$ of $R$, let
$R_{\cal{P}}$ be the localization of $R$ at ${\cal{P}}$. For an
open set $U\subseteq {\rm Spec}(R)$, we define ${\cal{O}}_{{\rm
Spec}(R)}(U)$ to be the set of functions $s: U\rightarrow
\coprod_{{\cal{P}}\in U}R_{\cal{P}}$, such that $s({\cal{P}})\in
R_{\cal{P}}$, for each ${\cal{P}}\in U$, and such that $s$ is
locally a quotient of elements of $R$: to be precise, we require
that for each ${\cal{P}}\in U$, there is a neighborhood $V$ of
${\cal{P}}$, contained in $U$, and elements $a,f\in R$, such that
for each ${\cal{Q}}\in  V$, $f\not\in {\cal{Q}}$, and
$s({\cal{Q}})=\frac{a}{f}$ in $R_{\cal{Q}}$ (see for example
Hartshorne \cite{Hartshorne},
 for definition and  basic properties of the sheaf
${\mathcal O}_{{\rm Spec}(R)}$).

Let $R$ be a ring. Then an $R$-module $M\neq 0$ is called a {\it
prime module} if  $rm=0$ for $m\in M$, $r\in R$ implies that $m=0$
or $rM=(0)$ (i.e., $r\in {\rm Ann}(M)$). We call a proper
submodule $P$ of an $R$-module $M$ to be a {\it prime submodule}
of $M$ if $M/P$ is a prime module, i.e., whenever $rm\in P$, then
either $m\in P$ or $rM\subseteq P$ for every $r\in  R$, $m\in  M$.
Thus $P$ is a prime submodule (or ${\cal{P}}$-prime submodule) of
$M$ if and only if ${\cal{P}}={\rm Ann}(M/P)$ is a prime ideal of
R and $M/P$ is a torsion free $R/{\cal{P}}$-module. This notion of
prime submodule was first introduced and systematically studied in
\cite{Dauns,Feller} and recently has received a good deal of
attention from several authors; see for examples
\cite{Behboodi1,Behboodi3,Behboodi4,Lu1,Lu5,Man,McCasland1,McCasland3}.
 The set of all prime submodules of $M$ is
called the {\it {\rm Spec}trum} of $M$ and denoted by ${\rm
Spec}(M)$. For any submodule $N$ of $M$, we have a set
$V(N)=\{P\in {\rm Spec}(M)~|~(N:M)\subseteq (P:M)\}$. Then the
sets $V(N)$, where $N$ is a submodule of $M$ satisfy the axioms
for the closed sets of a topology on {\rm Spec}$(M)$, called the
Zariski topology of $M$ (see for examples
\cite{Lu2,Lu3,McCasland2,Tekir}).

For an $R$-module $M$ let $\psi: {\rm Spec}(M)\rightarrow {\rm
Spec}(R/{\rm Ann}(M))$ defined by $\psi(P)=(P:M)/{\rm Ann}(M)$ for
every $P\in {\rm Spec}(M)$. $\psi$ is called the natural map of
${\rm Spec}(M)$. An $R$-module $M$ is called {\it primeful}  if
either $M=(0)$ or $M\neq (0)$ and the natural map $\psi: {\rm
Spec}(M)\rightarrow {\rm Spec}(R/{\rm Ann}(M))$ is surjective.
This notion of primeful module has been extensively studied by Lu
\cite{Lu4}).

Let $M$ be a primeful faithful $R$-module.   In \cite{Tekir}, the
author obtained  an $R$-module $\mathcal{O}_X(U)$ for each open
set $U$ in $X={\rm Spec}(M)$ such that  $\mathcal{O}_X$ is a sheaf
of modules over $X$. In fact, $\mathcal{O}_X$ is a generalization
of the structure sheaf of rings to primeful faithful modules.

 The purpose of this paper and its sequel, is to develop the structure sheaf of
 rings to a wider class of modules  called {\it
$\mathbb{P}$-radical modules} ((modules $M$ satisfying the prime
radical condition ``$(\sqrt[p]{{\cal{P}}M}:M)={\cal{P}}$'' for
every prime  ideal ${\cal{P}}\supseteq {\rm Ann}(M)$, where
$\sqrt[p]{{\cal{P}}M}$ is the intersection of all prime submodules
of $M$ containing ${\cal{P}}M$ but, if $M$ has no prime submodule,
then
 $\sqrt[p]{{\cal{P}}M}=M$). In Section 1, we  introduce and  study  $\mathbb{P}$-radical
modules. Several characterizations of $\mathbb{P}$-radical are
given in Proposition 2.1.  It is shown that the  class of
$\mathbb{P}$-radical  modules contains the class of primeful
$R$-modules properly (Proposition 2.3 and Example 2.4). This
yields over any ring all free modules, all finitely generated
modules and all homogeneous semisimple modules lie in the class of
$\mathbb{P}$-radical modules. Moreover, if $R$ is a Noetherian
ring (or an integral domain), then every projective $R$-modules is
$\mathbb{P}$-radical (see Theorem 2.5 and Corollary 2.4(iii)). In
particular, if $R$ is an Artinian ring, then
 all $R$-modules
are $\mathbb{P}$-radical (Theorem 2.13). Although the converse  is
not true in general (Example 2.14), but we show that the converse
is true for certain classes of rings like  Noetherian rings and
integral domains (see Theorem 2.16 and Corollary 2.17). Also, an
$R$-module $M$ is called {\it $\mathbb{M}$-radical module} if
$(\sqrt[p]{{\cal{M}}M}:M)={\cal{M}}$; for every maximal ideal
${\cal{M}}$ containing ${\rm Ann}(M)$. Thus, we have  the
following
chart of implications for $M$:\vspace{1mm}\\
 \indent\indent {\it $M$ is finitely generated $\Rightarrow$ $M$ is  primeful $\Rightarrow$  $M$ is $\mathbb{P}$-radical $\Rightarrow$
$M$ is  $\mathbb{M}$-radical}\vspace{1mm}\\ We show that the two
concepts $\mathbb{P}$-radical and $\mathbb{M}$-radical are
equivalent for all $R$-modules if and only if $R$ is a Hilbert
ring (Theorem 2.11). Also the two concepts primful and
$\mathbb{M}$-radical are equivalent for all $R$-modules if and
only if $dim(R)=0$ (Theorem 2.12). Recall that an $R$-module $M$
is called a multiplication module if for every submodule $N$ of
$M$ there exists an ideal $I$ of $R$ such that $N=IM$ (see for
example \cite{Barnard, El-Bast} for more details). Proposition
2.18, suggests that for a  multiplication module $M$ the four
concepts ``finitely generated", ``primeful",
``${\mathbb{P}}$-radical" and ``${\mathbb{M}}$-radical"  are
equivalent. Also we proof  the Proposition 2.19  that is analogue
of Nakayama's Lemma for $\mathbb{P}$-radical modules. In Section 3
 semisimple primful modules, semisimple $\mathbb{P}$-Radical modules and  semisimple $\mathbb{M}$-Radical modules,  are fully investigated.
 For instance, it is shown that a semisimple $R$-module $M$  is a
 ${\mathbb{M}}$-radical if and only if there exists a submodule $N$ of $M$ such that $N\cong\bigoplus_{{\rm Ann}(M)\subseteq{\cal{M}}\in
Max(R)}R/{\cal{M}}$) (Proposition 3.5.)  Also, a semisimple
$R$-module $M$  is a
 ${\mathbb{P}}$-radical if and only if $M$ is a ${\mathbb{M}}$-radical module and $R/{\rm Ann}(M)$ is a Hilbert ring (Proposition 3.6).
  In Part II we shall continue the study of this
construction and, as an application, we show that the sheaf theory
of  spectrum of $\mathbb{P}$-radical modules (with Zariski
topology) resembles to that of rings.

\section{$\mathbb{P}$-Radical  and $\mathbb{M}$-Radical Modules }

  Unlike the rings with identity, not every $R$-module
contains a prime submodule; for example ${\Bbb{Z}}_{p^\infty}$ as
$\Bbb{Z}$-module does not contain a prime submodule (see [5] or
[11]). Let $N$ be a proper submodule of an $R$-module $M$. Then
the {\it prime radical} $\sqrt[p]{N}$  is the intersection of all
prime submodules of $M$ containing $N$ or, in case there are no
such prime submodules, $\sqrt[p]{N}$  is $M$. Clearly
$V(N)=V(\sqrt[p]{N})$. We note that, for each ideal $I$ of $R$,
$\sqrt[p]{I}=\sqrt{I}$ (the intersection of all prime ideals of
$R$ containing $I$).  The prime radical of submodules are studied
by several authors; see for example
\cite{Behboodi2,Behboodi3,Lu5}.

The problem asking which $R$-module $M$ satisfies the equality
$\sqrt[p]{IM}=\sqrt IM$ for every ideal $I$ containing ${\rm
Ann}(M)$, had been investigated for finitely generated modules in
\cite{Lu5}. Also, in \cite{Lu4}, the author, extend the
investigation to  primeful flat content modules (e.g. free
modules), and primeful flat modules over rings with Noetherian
spectrum.

In this article we introduce a slightly different of the above
equality; that is the equality $(\sqrt[p]{IM}:M)=\sqrt I$ for
every ideal $I$ containing ${\rm Ann}(M)$. This prime radical
condition of modules play a key role in our investigation for give
a  sheaf construction for modules.

 The following  proposition offer several characterizations
of $R$-modules  $M$ satisfies the equality $(\sqrt[p]{IM}:M)=\sqrt
I$ for every ideal $I$ containing ${\rm Ann}(M)$.
\begin{Pro} For an $R$-module $M$ the following statements are equivalent:\vspace{3mm}\\
{\rm (1)} {\it $(\sqrt[p]{IM}:M)=\sqrt{I}$  for each  ideal $I\supseteq Ann(M)$.}\vspace{1mm}\\
{\rm(2)} {\it  $(\sqrt[p]{{\mathcal{P}}M}:M)={\mathcal{P}}$
 for each  prime ideal $\mathcal{P}\supseteq Ann(M)$.}\vspace{1mm}\\
{\rm(3)} {\it $\sqrt{I}=\bigcap_{P\in V(IM)}(P:M)$, for every
ideal
$I\supseteq Ann(M)$.}\vspace{1mm}\\
{\rm (4)} {\it ${\mathcal{P}}=\bigcap_{P\in
V({\mathcal{P}}M)}(P:M)$, for every prime ideal
${\mathcal{P}}\supseteq Ann(M)$.}
\end{Pro}
\noindent {Proof.} $(1)\Rightarrow(2)$ is trivial.\\
$(2)\Rightarrow(1)$. Let $I\supseteq {\rm Ann}(M)$. Clearly,
$I\subseteq (IM:M)\subseteq (\sqrt[p]{IM}:M)$ and
$\sqrt[p]{IM}:M)$ is a radical ideal (i.e., is an intersection of
all prime ideals). Thus $\sqrt{I}\subseteq (\sqrt[p]{IM}:M)$. On
the other hand,
$$(\sqrt[p]{IM}:M)\subseteq \bigcap_{{\mathcal{P}}\in
V(I)}(\sqrt[p]{{\mathcal{P}}M}:M)=\bigcap_{{\mathcal{P}}\in
V(I)}{\mathcal{P}}=\sqrt{I}.$$ Thus,
$(\sqrt[p]{IM}:M)=\sqrt{I}$.\\
 $(1)\Rightarrow(3)$. Let $I\supseteq {\rm Ann}(M)$.
It is easy to check
that $$\sqrt{I}=(\sqrt[p]{IM}:M)=((\bigcap_{P\in V(IM)}P):M)\\
=\bigcap_{P\in V(IM)}(P:M).$$\\
$(3)\Rightarrow(4)$ is trivial.\\
$(4)\Rightarrow(2)$.  Suppose that ${\mathcal{P}}$ is a prime
ideal of $R$  such that ${\mathcal{P}}\supseteq {\rm Ann}(M)$.
Since $\bigcap_{P\in V({\mathcal{P}}M)}P=\sqrt[p]{{\mathcal{P}}M}$
and $\bigcap_{P\in V({\mathcal{P}}M)}(P:M)=(\bigcap_{P\in
V({\mathcal{P}}M)}P:M)$, we conclude that
 ${\mathcal{P}}=(\sqrt[p]{{\mathcal{P}}M}:M)$. $\square$
\begin{Def}. {\rm  Let $M$ be  an $R$-module $M$. Then $M$ is called  a {\it $\mathbb{P}$-radical module}  if $M$ satisfies the equivalent
conditions listed in the above Proposition 2.1.}
\end{Def}
 The following  proposition shows that for any ring $R$, all primeful modules are $\mathbb{P}$-radical.
\begin{Pro} Let $R$ be a ring. Then any
primeful $R$-module $M$ is a $\mathbb{P}$-radical module.
\end{Pro}
\noindent {Proof.} Let $\mathcal{P}\supseteq {\rm Ann}(M)$. Since
$M$ is a primeful module, there exists $P\in {\rm Spec}(M)$ such
that $(P:M)=\mathcal{P}$. It follows that
$\sqrt[p]{\mathcal{P}M}\subseteq P$ and
$\mathcal{P}=\sqrt{\mathcal{P}}\subseteq
(\sqrt[p]{\mathcal{P}M}:M)$. Thus
$$\mathcal{P}=\sqrt{\mathcal{P}}\subseteq
(\sqrt[p]{\mathcal{P}M}:M)\subseteq (P:M)=\mathcal{P},$$
 and hence  $(\sqrt[p]{\mathcal{P}M}:M)=\mathcal{P}$. Thus    $M$ is a $\mathbb{P}$-radical module. $\square$\\

Let $R$ be a ring. It is easy to  see that every free $R$-module
$M$ is primeful. Also in [16, Theorem 2.2], it is shown that every
finitely generated $R$-module $M$ is primeful. In particular, if
$R$ is a domain, then every  projective $R$-module is primeful
(see [16, Corollary 4.3]). Thus by using these facts and above
proposition,  we have the following corollary.
\begin{Cor} Let $R$ be a ring.\vspace{2mm}\\
{\rm(i)} Every free  $R$-module is
$\mathbb{P}$-radical.\\
{\rm(ii)} Every finitely generated $R$-module is
$\mathbb{P}$-radical.\\
 {\rm(iii)}  If $R$ is a domain, then every projective $R$-module
is $\mathbb{P}$-radical.
\end{Cor}

Next, we show that Corollary 2.4 (iii) is also true when we
  replace ``$R$ is a domain" with ``$R$ is a Noetherian ring".
\begin{The} Let $R$ be a Noetherian ring. Then every projective $R$-module
is $\mathbb{P}$-radical.
\end{The}
\noindent {\bf Proof.} Suppose that $M$ is a projective $R$-module
and  ${\cal{P}}\supseteq {\rm Ann}(M)$. We claim that  that
${\cal{P}}M\neq M$. If ${\cal{P}}M=M$, then we put
${\cal{A}}=\{I\unrhd R~|~IM\neq M~{\rm and}~{\rm Ann}(M)\subseteq
I\subseteq {\cal{P}}\}$. Then ${\rm Ann}(M)\in {\cal{A}}$ and so
${\cal{A}}\neq \emptyset$. Since $R$ is Noetherian, ${\cal{A}}$
has a maximal element, say ${\cal{P}}_0$. Then ${\cal{P}}_0$ is a
prime ideal of $R$, otherwise, there exist   $a$, $b\in R\setminus
{\cal{P}}_0$ such that $ab\in {\cal{P}}_0$. It follows that
$({\cal{P}}_0+Ra)M=({\cal{P}}_0+Rb)M=M$ and so
$M=({\cal{P}}_0+Ra)({\cal{P}}_0+Rb)M\subseteq {\cal{P}}_0M$, a
contradiction. It is well-known  that for each projective
$R$-module $M$ and each ideal $I$ of $R$, the factor module $M/IM$
is also projective as an $R/I$-module. Thus
${\bar{M}}:=M/{\cal{P}}_0M$ is  projective as an
${\bar{R}}:=R/{\cal{P}}_0$-module. Since ${\bar{R}}$ is a domain,
by Corollary 2.4 (iii), ${\bar{M}}$ is a $\mathbb{P}$-radical
${\bar{R}}$-module. If $r\in R\setminus{\cal{P}}_0$, then
$(Rr+{\cal{P}}_0)M=M$ and it follows that  $r+{\cal{P}}_0\not\in
{\rm Ann}_{{\bar{R}}}({\bar{M}})$, i.e., ${\rm
Ann}_{{\bar{R}}}({\bar{M}})=(0)$.
 Since  ${\bar{\cal{P}}}:={\cal{P}}/{\cal{P}}_0$ is a prime ideal of ${\bar{R}}$ we
 must have
 $(\sqrt[p]{{\bar{\cal{P}}}{\bar{M}}}:{\bar{M}})={\bar{\cal{P}}}$.
 But, the equality   ${\cal{P}}M=M$ implies that ${\bar{\cal{P}}}{\bar{M}}={\bar{M}}$ and so
  $(\sqrt[p]{{\bar{\cal{P}}}{\bar{M}}}:{\bar{M}})=R$, a contradiction. Thus  ${\cal{P}}M\neq M$ and so ${\cal{P}}M$
 is a proper submodule of $M$.
   Suppose $F=M\oplus L$ were $F$
is a free $R$-module and $L$ is a submodule of $F$. Clearly
${\cal{P}}F$ is a prime submodule of $F$, i.e., $F/{\cal{P}}F$ is
a prime $R$-module and ${\rm Ann}(F/{\cal{P}}F)={\cal{P}}$. Since
$F/{\cal{P}}F\cong M/{\cal{P}}M\oplus L/{\cal{P}}L$, we conclude
that $M/{\cal{P}}M$ is also a prime $R$-module with ${\rm
Ann}(M/{\cal{P}}M)={\rm Ann}(F/{\cal{P}}F)={\cal{P}}$, i.e.,
${\cal{P}}M$ is  a prime submodule of $M$ with
$({\cal{P}}M:M)={\cal{P}}$. Thus
$(\sqrt[p]{{\cal{P}}M}:M)=({\cal{P}}M:M)={\cal{P}}$. Therefore $M$
is a $\mathbb{P}$-radical $R$-module. $\square$\\

 The following example shows  that the converse of Proposition 2.3
 is not  true in general (see also next Proposition  2.7).  Thus, the  class of $\mathbb{P}$-radical  modules  contains the class of primeful
$R$-modules properly.
\begin{Examp} {\rm (See also \cite{Lu4}, P. 136, Example 1.) Let
$M=\bigoplus_{p\in \Omega}\Bbb{Z}/p\Bbb{Z}$ as a $\Bbb{Z}$-module
where $\Omega$ is the set of prime integers. Clearly  ${\rm
Ann}(M)=0$ and according to \cite{Lu4}, for any non-zero prime
ideal $(p)$ of $\Bbb{Z}$, $pM$ is a $(p)$-prime submodule of $M$,
but it doesn't have any $(0)$-prime submodule. Thus $M$ is not a
primeful $\Bbb{Z}$-module. Clearly the zero submodule $(0)$ of $M$
is an intersection of maximal submodules of $M$  and hence
$\sqrt[p]{(0)}=(0)$. Thus $(\sqrt[p]{(0)M}:M)=((0):M)={\rm
Ann}(M)=(0)$. Also,  if $(q)$ is a nonzero prime (maximal) ideal
of $\Bbb{Z}$, then $(q)M=\bigoplus_{q\neq p\in
\Omega}\Bbb{Z}/p\Bbb{Z}\neq M$. It follows that
$(q)=(\sqrt[p]{(q)M}:M)$. Thus   $M$ is a $\mathbb{P}$-radical
module and it is not a primeful module.}
\end{Examp}
\begin{Pro} Let $R$ be a    ring. Then there exists an $R$-module $M$ which
is a $\mathbb{P}$-radical module  but it is not a primeful module
if and only if there exist prime ideals $\mathcal{P}$ and
${\{{\mathcal{P}}_i\}}_{i\in I}$
  of $R$ such that  $\mathcal{P}\subsetneqq {\mathcal{P}}_i$ (for each $i\in
  I$)
  and $\mathcal{P}=\bigcap_{i\in I}{\mathcal{P}}_i$.
  \end{Pro}
\noindent {\bf Proof.} $(\Rightarrow)$ Let  $M$ be an  $R$-module
which is a $\mathbb{P}$-radical module  but it is not a primeful
module. Assume that  $\mathcal{P}$ is a prime ideal of $R$ such
that $M$ has not any $\mathcal{P}$-prime submodule. Thus for any
prime submodule $N$ of $M$ with $\mathcal{P}M\subseteq N$, we must
have $\mathcal{P}\subsetneqq {\mathcal{P}}_N:=(N:M)$. Since $M$ is
a $\mathbb{P}$-radical module,
$$\mathcal{P}=(\sqrt[p]{\mathcal{P}M}:M)=(\bigcap_{\mathcal{P}M\subseteq
N\in {\rm Spec}(M)}N:M)=\bigcap_{\mathcal{P}M\subseteq N\in {\rm
Spec}(M)}(N:M)=\bigcap_{\mathcal{P}M\subseteq N\in {\rm
Spec}(M)}{\mathcal{P}}_N.$$
 $(\Leftarrow)$ Without loss of generality we can assume that
  $\mathcal{P}=\bigcap_{i\in I}{\mathcal{P}}_i$ and for each prime ideal $\mathcal{Q}\supsetneqq {\mathcal{P}}$ of $R$ there exists $i\in I$ such that
  $\mathcal{Q}={\mathcal{P}}_i$.
 Let $M=\bigoplus_{i\in I}R/{\mathcal{P}}_i$ as
$R$-module.  Clearly ${\rm Ann}(M)=\bigcap_{i\in
I}{\mathcal{P}}_i=\mathcal{P}$ and  for each $j\in I$,
$N_j=\bigoplus_{j\neq i\in I}R/{\mathcal{P}}_i$ is a
${\mathcal{P}}_j$-prime submodule of $M$  with
${\mathcal{P}M}\subseteq {\mathcal{P}}_jM\subseteq N_j$.  Thus for
each prime ideal ${\mathcal{P}}_j\supsetneqq {\rm
Ann}(M)={\mathcal{P}}$,
$(\sqrt[p]{{\mathcal{P}}_jM}:M)={\mathcal{P}}_j$ (since
${\mathcal{P}}_j\subseteq (\sqrt[p]{{\mathcal{P}}_jM}:M)\subseteq
(N_j:M)={\mathcal{P}}_j$).  On the other hand
$$(\sqrt[p]{{\mathcal{P}}M}:M)\subseteq (\bigcap_{i\in I}N_i:M)=\bigcap_{i\in I}(N_i:M)=\bigcap_{i\in I}{\mathcal{P}}_i=\mathcal{P}.$$
Thus $M$ is a $\mathbb{P}$-radical module.  We  claim that $M$
doesn't have any ${\mathcal{P}}$-prime submodule, for if not, let
$N$ be a ${\mathcal{P}}$-prime submodule of $M$ with
$(N:M)={\mathcal{P}}$. Since $N\neq M$, there exists $j\in I$ such
that $(\cdots,0,1+{\mathcal{P}}_j,0,\cdots)\not\in N$. Since
${\mathcal{P}}_j(\cdots,0,1+{\mathcal{P}}_j,0,\cdots)\in N$, we
must have  ${\mathcal{P}}_jM\subseteq N$,  i.e.,
${\mathcal{P}}_j\subseteq {\mathcal{P}}$, a contradiction.  Thus
 $M$  is
not a  primeful $R$-module.$~\Box$\\

By Proposition 2.1 and Definition 2.2, an $R$-module $M$ is a
 $\mathbb{P}$-radical module if and only if $(\sqrt[p]{{\mathcal{P}}M}:M)={\mathcal{P}}$
 for each  prime ideal $\mathcal{P}\supseteq {\rm Ann}(M)$. Now we have to adapt the notion of $\mathbb{M}$-radical modules which generalized $\mathbb{P}$-radical modules.

\begin{Def}. {\rm  Let $M$ be  an $R$-module $M$. Then $M$ is called  an {\it $\mathbb{M}$-radical
module} (or a {\it Maxful module}) if
$(\sqrt[p]{{\mathcal{M}}M}:M)={\mathcal{M}}$
 for each  maximal  ideal $\mathcal{M}\supseteq {\rm Ann}(M)$.}
\end{Def}

 The following evident proposition offers several
other characterizations of $\mathbb{M}$-radical modules.

\begin{Pro} Let $R$ be a  ring and $M$ a nonzero $R$-module. Then the following statements  are equivalent:\vspace{3mm}\\
{\rm (1)}  $M$ is an $\mathbb{M}$-radical module.\vspace{1mm}\\
{\rm (2)}  ${\cal{M}}M\neq M$ for every maximal ideal
${\cal{M}}\supseteq Ann(M)$.\vspace{1mm}\\
{\rm (3)}  ${\cal{P}}M\neq M$ for every prime  ideal
${\cal{P}}\supseteq Ann(M)$. \vspace{1mm}\\
 {\rm (4)}  For every
maximal ideal ${\mathcal{M}}\supseteq Ann(M)$ there exists a
maximal  submodule $P$ of $M$ \indent such that
$(P:M)={\mathcal{M}}$.\vspace{1mm}\\
{\rm(5)}  For every maximal   ideal ${\mathcal{M}}\supseteq
Ann(M)$ there exists a prime  submodule $P$ of $M$ such \indent
that $(P:M)={\mathcal{M}}$.
\end{Pro}
 \noindent {\bf Proof.} $(1)\Rightarrow(2)$. Clearly for each  maximal  ideal $\mathcal{M}\supseteq {\rm Ann}(M)$,
 the equality
  $(\sqrt[p]{{\mathcal{M}}M}:M)={\mathcal{M}}$ implies that ${\cal{M}}M\neq M$. \\
  $(2)\Rightarrow(1)$.  Suppose  ${\cal{M}}M\neq M$ for every maximal ideal
${\cal{M}}\supseteq {\rm Ann}(M)$. Then for each maximal ideal
${\cal{M}}\supseteq {\rm Ann}(M)$, ${\cal{M}}M$ is a prime
submodule of $M$ and so $\sqrt[p]{{\mathcal{M}}M}={\mathcal{M}}M$.
It follows that  ${\mathcal{M}}\subseteq
 (\sqrt[p]{{\mathcal{M}}M}:M)=({{\mathcal{M}}M}:M)={\mathcal{M}}$, i.e.,  $(\sqrt[p]{{\mathcal{M}}M}:M)={\mathcal{M}}$.  Thus $M$ is an $\mathbb{M}$-radical module.\\
$(2)\Leftrightarrow(3)$ is clear.\\
$(2)\Rightarrow(5)$. Suppose that ${\mathcal{M}}\supseteq {\rm
Ann}(M)$ is a maximal ideal. Then ${\cal{M}}M\neq M$ implies that
${\cal{M}}M$ is a prime submodule with
$({\cal{M}}M:M)={\cal{M}}$.\\
$(5)\Rightarrow(2)$ is clear.\\
 $(4)\Rightarrow(5)$ is clear
(since every
maximal submodule is a prime submodule).\\
 $(5)\Rightarrow(4)$. Suppose that $(P:M)={\mathcal{M}}$ where $P$ is a prime submodule of $M$ and  ${\mathcal{M}}\supseteq
{\rm Ann}(M)$ is a maximal ideal. Then $M/P$ is an
$R/{\mathcal{M}}$-vector space and hence $M/P$ has a maximal
$R/{\mathcal{M}}$-subspace say $K/P$. Clearly $K\subseteq M$ is a
maximal $R$-submodule and $(K:M)={\mathcal{M}}$. $\square$

A commutative ring $R$ is called a {\it Hilbert ring}, also {\it
Jacobson} or {\it  Jacobson-Hilbert} ring, if  every prime ideal
of  $R$ is the intersection of maximal ideals. The class of
commutative Hilbert rings is closed under forming finite
polynomial rings. For an integral extension $R\subseteq S$, the
ring $S$ is Hilbert if and only if $R$ is a Hilbert ring. The main
interest in Hilbert rings in commutative algebra and algebraic
geometry is their relation with Hilbert's Nullstellensatz (see
Goldman \cite{Goldman},  Cortzen and Small \cite{Cortzen}, Theorem
1,  and \cite{Eisenbud}, Theorem 4.19, for more ditals).

Next, we will to show that the two concepts $\mathbb{P}$-radical
and $\mathbb{M}$-radical are equivalent for all $R$-modules if and
only if $R$ is a Hilbert ring. As we have observed earlier some
modules $M$ have no prime submodules (for example
${\Bbb{Z}}_{p^\infty}$ as $\Bbb{Z}$-module) and we call them {\it
primeless}. We recall that a module $M$ over a domain $R$ is
called {\it divisible} if $rM=M$ for each $0\neq r\in R$ and $M$
is called {\it torsion} if ${\rm Ann}(m)\neq 0$ for each $m\in M$.

 We need the following evident lemma.

\begin{Lem}  Let $R$ be a   domain. Then any torsion divisible $R$-module is
primeless.
\end{Lem}
\begin{The} Let $R$ be a
ring. Then the two concepts $\mathbb{P}$-radical and
$\mathbb{M}$-radical are equivalent for all $R$-modules if and
only if $R$ is a Hilbert ring.
\end{The}
\noindent {\bf Proof.} $(\Rightarrow)$. Let  every
$\mathbb{M}$-radical module is a $\mathbb{P}$-radical module. To
obtain a contradiction, suppose that $\mathcal{P}$ is a prime
ideal of $R$ such that it is not an intersection of maximal
ideals, i.e., $\mathcal{P}\neq
\bigcap_{\mathcal{P}\subseteq\mathcal{M}\in Max(R)}\mathcal{M}$.
It is easy to check that every  $\mathbb{M}$-radical
$R/\mathcal{P}$-module is also a $\mathbb{P}$-radical
$R/\mathcal{P}$-module. Let
 $S=R/\mathcal{P}$ and   $Q$ be the field of fraction of $S$.
 Since $\mathcal{P}$ is not a maximal ideal, $S$ is not a field
 and so $S\neq Q$.
 Suppose $K$ is a nonzero proper $S$-submodule of $Q$. Then  $L=Q/K$ is a torsion divisible $S$-module and so by Lemma 2.10, $L$
 is a primeless $S$-module. We put
$$M=(\bigoplus_{\mathcal{M}\in
Max(S)}S/\mathcal{M})\oplus L$$ as an $S$-module. Clearly $M$ is a
$\mathbb{M}$-radical $S$-module, since   for each
${\mathcal{M}}_1\in Max(S)$,
$${\mathcal{M}}_1M=(\bigoplus_{{\mathcal{M}}_1\neq\mathcal{M}\in
Max(S)}S/\mathcal{M})\oplus L$$ is a prime $S$-module with
$({\mathcal{M}}_1M:M)={\mathcal{M}}_1$. On the other hand we claim
that every prime $S$-submodule of $M$ is of the above form.   To
see this, let $P$ be a prime  $S$-submodule of $M$. We claim that
$\bigoplus_{\mathcal{M}\in Max(S)}S/\mathcal{M}\not\subseteq P$,
otherwise we must have $M/P\cong L/T$ for some proper submodule
$T$ of $L$, a contradiction (since   $L$
 is a primeless $S$-module).  Thus there exists ${\mathcal{M}}_1\in Max(S)$
such that $(0,\cdots,0,1+{\mathcal{M}}_1,0\cdots)\not\in P$. Since
 ${\mathcal{M}}_1(0,\cdots,0,1+{\mathcal{M}}_1,0\cdots)\subseteq  P$,  ${\mathcal{M}}_1M\subseteq N$, and hence ${\mathcal{M}}_1M=N$.
 It follows that $$\sqrt[p]{(0)}=\bigcap_{\mathcal{M}\in
Max(R)}\mathcal{M}M=L.$$  Since ${\rm Ann}(_SL)=(0)$, we conclude
that ${\rm Ann}(_SM)=(0)$. Thus
$$(\sqrt[p]{(0)M}:M)=(\sqrt[p]{(0)}:M)=(L:M)=((0):\bigoplus_{\mathcal{M}\in
Max(S)}S/\mathcal{M})=\bigcap_{\mathcal{M}\in
Max(R)}{\mathcal{M}}.$$ But    $\bigcap_{\mathcal{M}\in
Max(R)}{\mathcal{M}}\neq (0)$ since  $\mathcal{P}\neq
\bigcap_{\mathcal{P}\subseteq\mathcal{M}\in Max(R)}\mathcal{M}$.
Thus $(\sqrt[p]{(0)M}:M)\neq (0)$, which is impossible (since $M$
is a $\mathbb{P}$-radical module). Thus $R$ is a Hilbert ring.\\
$(\Leftarrow)$. Since every $\mathbb{P}$-radical module is a
$\mathbb{M}$-radical module, it suffices to show that every
$\mathbb{M}$-radical module is a $\mathbb{P}$-radical module.
Assume that $M$ is  an $\mathbb{M}$-radical module and
$\mathcal{P}$ is a prime ideal of $R$ with $\mathcal{P}\supseteq
{\rm Ann}(M)$. Since $R$ is a Hilbert ring,
$\mathcal{P}=\bigcap_{\mathcal{P}\subseteq \mathcal{M}\in
Max(R)}\mathcal{M}$. On the other hand, since $M$ is an
$\mathbb{M}$-radical module,  $(\mathcal{M}M:M)=\mathcal{M}$ for
each ${\mathcal{P}}\subseteq \mathcal{M}\in Max(R)$.
 It follows that $$(\sqrt[p]{\mathcal{P}M}:M)\subseteq(\bigcap_{\mathcal{P}\subseteq\mathcal{M}\in Max(R)}{\mathcal{M}}M:M)=
 \bigcap_{\mathcal{P}\subseteq\mathcal{M}\in Max(R)}({\mathcal{M}}M:M)=\bigcap_{\mathcal{P}\subseteq\mathcal{M}\in Max(R)}{\mathcal{M}}=\mathcal{P}.$$
Since $\mathcal{P}\subseteq (\mathcal{P}M:M)\subseteq
(\sqrt[p]{\mathcal{P}M}:M)$, we conclude that
$(\sqrt[p]{\mathcal{P}M}:M)=\mathcal{P}$. Thus $M$  a
$\mathbb{P}$-radical module. $\square$\\

Let $M$ be an $R$-module. By  Proposition 2.3 and Proposition 2.9,
we have  the following chart of implications for $M$:

{\it $M$ is a primeful module $\Rightarrow$  $M$ is a
$\mathbb{P}$-radical module $\Rightarrow$   $M$ is a
$\mathbb{M}$-radical module}

In general, none of the implications is reversible. However, for
zero dimensional  rings, the  implications can be replaced by
equivalences, as we shall now show.

\begin{The} Let $R$ be a   ring.  Then the two concepts primful and
$\mathbb{M}$-radical are equivalent for all $R$-modules if and
only if $dim(R)=0$.
\end{The}
\noindent {\bf Proof.} Assume that the two concepts primful and
$\mathbb{M}$-radical are equivalent for all $R$-modules. Then by
Theorem 2.11, $R$ is a Hilbert ring. Suppose, contrary to our
claim, that $dim(R)\geq 2$ and ${\cal{P}}$ is a non-maximal  prime
ideal of $R$. Thus ${\cal{P}}=\bigcap_{{\cal{P}}\subseteq
{\cal{M}}\in Max(R)}{\cal{M}}$. Therefore,  by  Proposition 2.7,
there exists an $R$-module $M$ which is a $\mathbb{P}$-radical
($\mathbb{M}$-radical) module  but it is not a primeful module, a
contradiction. Thus $dim(R)=0$. The converse is clear. $\square$

The following result shows over an    Artinian ring, every module
is $\mathbb{P}$-radical (consequently, every module is primeful
and also every module is $\mathbb{M}$-radical).

\begin{The}
Let $R$ be an   Artinian ring. Then every $R$-module is
 a $\mathbb{P}$-radical module.
\end{The}
\noindent {\bf Proof.} Since $dim(R)=0$, by  Theorem 2.12,  the
three concepts primful, $\mathbb{P}$-radical and
$\mathbb{M}$-radical are equivalent for all $R$-modules. Since $R$
is  Artinian,
  $R=R_1\times\cdots\times R_n$, where
 $n\in\Bbb{N}$  and
each $R_i$ is an   Artinian local ring. First  give the proof  for
the case $n=1$, i.e., $R$ is a local ring with maximal ideal
$\cal{M}$. Suppose $M$ is a nonzero $R$-module. Then ${\rm
Ann}(M)\neq R$ and so ${\rm Ann}(M)\subseteq \cal{M}$. By
Proposition 2.9, it is suffices to show that ${\cal{M}}M\neq M$.
If $R$ is a domain (field), then ${\cal{M}}=(0)$ and so  the proof
is complete. Hence, to obtain a contradiction, suppose that $R$ is
not a domain and also ${\cal{M}}M=M$. The there are nonzero
elements $a$, $b\in R$ such that $ab=0$. Thus  $(Ra)(Rb)M=(0)$ and
so either $RaM\neq M$ or $RbM\neq M$. It follows that
${\cal{A}}:=\{I\unrhd R~|~IM\neq M~and~I\neq (0)\}$  is a nonempty
set of ideals of $R$. Since $R$ is Noetherian, ${\cal{A}}$ has a
maximal element, say $\cal{P}$. Since ${\cal{P}}M\neq M$,
$\cal{P}$ is not a prime (maximal) ideal of $R$ (i.e.,
${\cal{P}}\neq {\cal{M}}$). Thus there exist ideals $A$, $B$ of
$R$ such that ${\cal{P}}\subsetneqq A$, ${\cal{P}}\subsetneqq B$
and $AB\subseteq {\cal{P}}$. Thus the equalities  $AM=BM=M$
implies that $M=ABM\subseteq {\cal{P}}M$, a contradiction.
Therefore,  in the case $R$ is a local ring every $R$-module is
 a $\mathbb{P}$-radical module. Now assume  $n\geq 2$ and for each $i$ ($1\leq i\leq n$) let ${\cal{M}}_i$ be  the maximal ideal of the local ring $R_i$.
 Suppose $M$ is a nonzero $R$-module and  ${\rm Ann}(M)\subseteq \cal{M}$ where $\cal{M}$ is a maximal ideal of
 $R$.  Clearly $\cal{M}$ is the form $R_1\times\cdots\times R_{i-1}\times {\cal{M}}_i\times R_{i+1}\cdots\times R_n$ for some
 $i$. Without loss of generality we can assume that $i=1$, i.e., ${\cal{M}}={\cal{M}}_1\times R_2\times\cdots\times R_n$.
 Again,  by Proposition
2.9, it is suffices to show that ${\cal{M}}M=({\cal{M}}_1\times
R_2\times\cdots\times R_n)M\neq M$. On the contrary, suppose that
 $({\cal{M}}_1\times R_2\times\cdots\times R_n)M=M$. Let  $I=R_1\times (0)\times \cdots\times (0)$ and $J=(0)\times R_2\times \cdots\times
 R_n$. Then $I$, $J$  are ideals of $R$ with $J={\rm Ann}(I)$. Thus $R_1\cong R/J$ and  so ${\bar{M}}=IM$ is an
 unitary $R_1$-module (in fact, ${\bar{M}}=(R_1\times 0\times \cdots\times 0)M$ is an
 unitary $R_1$-module with $r_1{\bar{m}}$ defined
  to be $r_1(1,0,\cdots,0){\bar{m}}$  for  $r_1\in R_1$ and ${\bar{m}}\in \bar{M}$). We claim that ${\bar{M}}\neq (0)$, otherwise,
   $R_1\times (0)\times \cdots\times (0)\subseteq {\rm Ann}(M)\subseteq {\cal{M}}_1\times R_2\times\cdots\times R_n$, a contradiction.
   Thus ${\bar{M}}$ is a nonzero $R_1$-module and so by using  case $n=1$, we have
   ${\cal{M}}_1{\bar{M}}\neq{\bar{M}}$, i.e., $({\cal{M}}_1\times 0\times \cdots\times 0)M\neq (R_1\times 0\times \cdots\times
   0)M$. On the other hand,  for each $m\in M$,  $(1,0,\cdots,0)m\in M=({\cal{M}}_1\times R_2\times\cdots\times R_n)M$. Thus for each  $m\in M$,
    $$(1,0,\cdots,0)m=\sum_{j=1}^k(p_{1j},r_{2j},\cdots,r_{nj})m_j$$
    for some $k\in\Bbb{N}$, $m_j\in M$,  $p_{1j}\in {\cal{M}}_1$, $r_{ij}\in R_i$ were  $2\leq i\leq n$ and  $1\leq j\leq
    k$. Thus  $$(1,0,\cdots,0)m=(1,0,\cdots,0)^2m=\sum_{j=1}^k(p_{1j},0,\cdots,0)m_j\in ({\cal{M}}_1\times (0)\times\cdots\times (0))M,$$
 for each $m\in M$. It follows that $(R_1\times 0\times\cdots\times 0)M\subseteq ({\cal{M}}_1\times (0)\times\cdots\times (0))M$, i.e.,
  ${\cal{M}}_1{\bar{M}}={\bar{M}}$, a contradiction. $\square$ \\

The following example shows that in general the converse of
Theorem 2.13 is not true.

\begin{Examp}
 {\rm Let $K$ be a field, $D:=K[\{x_i : i\in\Bbb{N}\}]$ (a unique factorization domain) and
$R=K[\{x_i : i\in \Bbb{N}\}]/(\{x_ix_j : i,j\in
 \Bbb{N}\})$ where $(\{x_ix_j : i,j\in
 \Bbb{N}\})$ is the ideal of $D$ generated by the subset $\{x_ix_j : i,j\in
 \Bbb{N}\}\subseteq D$. Let $\bar{x}_k=x_k+(\{x_ix_j : i,j\in
 \Bbb{N}\})$ for each $k\in\Bbb{N}$ and  ${\cal{M}}=(\{\bar{x}_k :  k\in
 \Bbb{N}\})$. The ideal ${\cal{M}}$ is simply the image of the maximal ideal ${\cal{N}}=(\{x_k : k\in
 \Bbb{N}\})$ of $D$. Clearly ${\cal{M}}^2=(0)$ and so  $R$
is a local zero-dimensional ring, but is not Artinian
(Noetherian). Since ${\cal{M}}^2=(0)$,  for each nonzero
$R$-module $M$, ${\cal{M}}M\neq M$. Thus by Proposition 2.9 and
Theorem 2.12,  $M$ is a $\mathbb{P}$-radical module. Thus  every
$R$-module is $\mathbb{P}$-radical, but $R$ is not Artinian.}
\end{Examp}
A ring $R$ is called a  Max-ring (or a Bass ring) if every nonzero
$R$-module has a maximal submodule. Also, a ring $R$ is called a
P-ring if every nonzero $R$-module has a prime submodule. In
\cite{Behboodi4}, Theorem 3.9,  it is shown that commutative
P-rings coincide with Max-rings. Moreover, we have the following
lemma which is essentially Theorem 2 of \cite{Hamsher}.
\begin{Lem}  For a commutative ring $R$, the following conditions are
equivalent:\vspace{3mm}\\
{\rm (1)} $R$ is a max ring;\vspace{1mm}\\
 {\rm(2)} $R/J(R)$ is a regular ring and
$J(R)$ is a t-nilpotent ideal.
\end{Lem}
The following  theorem offer several characterizations of
Noetherian rings $R$ over which  every  module is
$\mathbb{P}$-radical.

\begin{The}  Consider the following statements for a   ring  $R$:\vspace{3mm}\\
{\rm (1)}  $R$ is an Artinian ring.\vspace{1mm}\\
{\rm (2)}  Every  $R$-module is  primeful.\vspace{1mm}\\
{\rm(3)}  Every $R$-module is $\mathbb{P}$-radical.\vspace{1mm}\\
{\rm(4)}  Every $R$-module is  $\mathbb{M}$-radical.\vspace{1mm}\\
{\rm(5)}  $R$ is a Max-ring.\vspace{1mm}\\
{\rm(6)}  $R$ is a P-ring.\vspace{1mm}\\
{\rm(7)}  $dim(R)=0$.\vspace{2mm}\\
Then
$(1)\Rightarrow(2)\Leftrightarrow(3)\Leftrightarrow(4)\Rightarrow(5)\Leftrightarrow(6)\Rightarrow(7)$.
Consequently, when $R$ is a Noetherian ring (or a domain), all the
six statements are equivalent.
\end{The}
\noindent {\bf Proof.} $(1)\Rightarrow(2)$ is  by Theorem 2.13.\\
$(2)\Rightarrow(3)$ is by Proposition 2.3. \\
$(3)\Rightarrow(4)$ is clear. \\
$(4)\Rightarrow(5)$. Assume that every $R$-module is
$\mathbb{M}$-radical. Let
 $M$ be  a nonzero $R$-module. Then ${\rm Ann}(M)\neq R$ and so there
exists a maximal ideal $\cal{M}$ of $R$ such that ${\rm
Ann}(M)\subseteq \cal{M}$. Since $M$ is maxful, there exists a
prime submodule $P$ of $M$ with $(P:M)=\cal{M}$. Thus $M/P$ is an
$R/\cal{M}$-module ($R/\cal{M}$-vector space), and hence, $M/P$
has a maximal $R/\cal{M}$-submodule, say $K/P$. It is easy to see
that
$K\lneqq M$ is a maximal $R$-submodule. Thus every nonzero $R$-module has a maximal submodule, i.e., $R$ is a Max-ring.  \\
$(5)\Leftrightarrow(6)$ is by \cite{Behboodi4}, Theorem 3.9.\\
$(6)\Rightarrow(7)$. Let $R$ be a P-ring. Suppose that $\cal{P}$
is a non-maximal prime ideal of $R$.  Set $R'=R/{\cal{P}}$  and
let $K$ be the field of fractions of $R'$. Since $R'$ is not a
field, $R'\neq K$ and $K$ is a divisible $R'$-module. It follows
that $K/R'$ is a nonzero torsion divisible $R'$-module. Thus by
Lemma 2.10, $K/R'$ is a primeless $R'$-module. Now it is easy to
see that $K/R'$ is a primeless $R$-module. Thus $R$
is not a P-ring, a contradiction.\\
$(4)\Rightarrow(2)$. Assuming (4) to hold. We conclude from
$(4)\Rightarrow(7)$ that $dim(R)=0$.  Theorem 2.12 now shows that
every  $R$-module is  primeful.\\
Finally, if $R$  is a Noetherian ring, then $dim(R)=0$ if and only
if $R$ is an Artinian ring. Thus $(6)\Rightarrow(1)$, when $R$  is
a Noetherian ring.  $\square$

The following is now immediate.

\begin{Cor}  Let $R$ be a domain. Then the following statements are equivalent:\vspace{3mm}\\
{\rm (1)}  Every  $R$-module is  primeful.\vspace{1mm}\\
{\rm(2)}  Every $R$-module is $\mathbb{P}$-radical.\vspace{1mm}\\
{\rm(3)}  Every $R$-module is  $\mathbb{M}$-radical.\vspace{1mm}\\
{\rm(4)}  $R$ is a field (i.e., $R$ is an Artinan domain).
\end{Cor}

The following proposition  suggests that for a  multiplication
module $M$ the four concepts ``finitely generated",
 ``primeful", ``${\mathbb{P}}$-radical" and  ``${\mathbb{M}}$-radical"  are
equivalent.

\begin{Pro} Consider the following statements for a nonzero
$R$-module $M$:\vspace{3mm}\\
{\rm  (1)}  $M$ is finitely generated.\vspace{1mm}\\
{\rm  (2)}  $M$ is primeful.\vspace{1mm}\\
{\rm  (3)} $M$ is a ${\mathbb{P}}$-radical module.\vspace{1mm}\\
{\rm  (4)}   $({\cal{P}}M : M)={\cal{P}}$ for every
prime ideal ${\cal{P}}\supseteq Ann(M)$.\vspace{1mm}\\
 {\rm  (5)} $M$ is a ${\mathbb{M}}$-radical module.\vspace{3mm}\\
 Then {\rm $(1)\Rightarrow (2)\Rightarrow (3)\Rightarrow (4)\Rightarrow(5)$}. When $M$ is a multiplication module,
 then {\rm $(5)\Rightarrow (1)$}.
\end{Pro}
\noindent {\bf Proof.} $(1)\Rightarrow(2)$ is  by \cite{Lu4},
Proposition 3.8. \\
$(2)\Rightarrow(3)$ is by Proposition 2.3. \\
$(3)\Rightarrow(4)$. Since ${\cal{P}}\subseteq ({\cal{P}}M :
M)\subseteq (\sqrt[p]{{\cal{P}}M} : M)$ for every prime ideal
${\cal{P}}\supseteq {\rm Ann}(M)$, the proof is clear.\\
$(4)\Rightarrow(5)$. For every
prime ideal ${\cal{P}}\supseteq {\rm Ann}(M)$ the equality $({\cal{P}}M : M)={\cal{P}}$ implies that ${\cal{P}}M\neq M$. Thus $M$ is a ${\mathbb{M}}$-radical module. \\
When $M$ is a multiplication module, then $(5)\Rightarrow (1)$ is
by \cite{Lu4}, Proposition 3.8. $\square$

We conclude this section  with the next analogue of Nakayama's
Lemma.
\begin{Pro}
Let $M$ be a ${\mathbb{M}}$-radical $R$-module. Then $M$ satisfies
the following condition (NAK): If $I$ is an ideal of $R$ contained
in the Jacobson radical $J(R)$ such that $IM=M$, then $M=(0)$.
\end{Pro}
\noindent {\bf Proof.} Suppose that $M\neq (0)$. Then ${\rm
Ann}(M)\neq R$. If ${\cal{M}}$ is any maximal ideal containing
${\rm Ann}(M)$, then $I\subseteq {\cal{M}}$ and $IM=M={\cal{M}}M$,
a contradiction. $\square$

\section{ Characterization of   Semisimple $\mathbb{P}$-Radical   Modules }

Recall that for a module $M$, the {\it socle} of $M$ (denoted by
soc$(M)$) is the sum of all simple (minimal) submodules of $M$. If
there are no minimal submodules in $M$ we put soc$(M)=(0)$. Thus
$M$ is a semisimple module if soc$(M)= M$. Also, a semisimple
module $M$ is called homogeneous if any two simple submodules of
$M$ are isomorphic. It is  easy to check that an $R$-module $M$ is
homogeneous semisimple if and only if ${\rm Ann}(M)$ is a maximal
ideal.

Next, we aim to characterize ${\mathbb{P}}$-radical semisimple
modules. First we need the following definition.

\begin{Def}. {\rm  Let $R$ be a ring. A  semisimple $R$-module $M$ is called  {\it
full semisimple}  if for each maximal ideal   ${\cal{M}}\supseteq
{\rm Ann}(M)$ the simple $R$-module $R/{\cal{M}}$ may be embedded
 in $M$ (i.e., there exists a submodule $N$ of $M$ such that $N\cong\bigoplus_{{\rm Ann}(M)\subseteq{\cal{M}}\in
Max(R)}R/{\cal{M}}$).}
\end{Def}

\begin{Examp} {\rm Let $M=\bigoplus_{p\in \Omega}\Bbb{Z}/p\Bbb{Z}$ as a $\Bbb{Z}$-module
where $\Omega$ is the set of prime integers. Clearly $M$ is a full
semisimple  $\Bbb{Z}$-module. But the semisimple $\Bbb{Z}$-module
$M_1=\bigoplus_{2\neq p\in \Omega}\Bbb{Z}/p\Bbb{Z}$ is not full
semisimple since ${\rm Ann}(M_1)=(0)\subseteq 2\Bbb{Z}$  and
$\Bbb{Z}/2\Bbb{Z}$ is not a submodule of $M_1$. }
\end{Examp}

The proof of the following result is straightforward and  left to
the reader.

\begin{Pro}  Let $M$ be a semisimple $R$-module such that $Ann(M)$ is a finite intersection of maximal ideals.
  Then $M$ is full semisimple.  In particular,  all  finitely generated semisimple modules as well as all  homogenous semisimple
  modules are full semisimple.
\end{Pro}

\begin{Lem}  Let $M$ be an $R$-module with nonzero socle.  Then $M$ is a prime module if and only if $M$ is a homogeneous
semisimple module (i.e., $Ann(M)$ is a maximal ideal).
\end{Lem}
\noindent {\bf Proof.}   Let $M$ be a prime module with nonzero
socle and let $Rm$ be a simple submodule of $M$. Then ${\rm
Ann}(m)={\rm Ann}(M)=P$ and hence $P$ is a maximal ideal of $R$.
Since ${\rm Ann}(m)={\rm Ann}(m')$ for each $0\neq m'\in M$, $M$
is a homogeneous
semisimple $R$-module. The converse is evident. $\square$\\

Now we are in position to show that the two   concepts
 $\mathbb{M}$-radical and full semisimple  are
equivalent for semisimple modules.

\begin{Pro}  Let $M$ be a semisimple $R$-module.  Then $M$ is  ${\mathbb{M}}$-radical  if and only if
 $M$ is full semisimple.
\end{Pro}
\noindent {\bf Proof.} Since $M$ is a semisimple $R$-module, we
can assume that $M=\bigoplus_{i\in I}R/{\cal{M}}_i$ where $I$ is
an index set and each  ${\cal{M}}_i$ is a maximal ideal of $R$.\\
 $(\Rightarrow)$. Clearly ${\rm Ann}(M)=\bigcap_{i\in I}{\cal{M}}_i$.
Suppose ${\cal{M}}\supseteq {\rm Ann}(M)$ is a maximal ideal of
$R$. If ${\cal{M}}\neq  {\cal{M}}_i$ for each $i\in I$, then
${\cal{M}}(R/{\cal{M}}_i)=R/{\cal{M}}_i$ for each $i\in I$. It
follows that ${\cal{M}}M=M$, contrary to (2) (see Proposition
2.9). Thus ${\cal{M}}={\cal{M}}_i$ for some  $i\in I$ and hence
$R/{\cal{M}}$ may be embedded  in $M$.\\
$(\Leftarrow)$. Suppose ${\cal{M}}\supseteq {\rm Ann}(M)$ is a
maximal ideal of $R$. Since $M$ is full semisimple,
${\cal{M}}={\cal{M}}_i$ for each $i\in I$ and so ${\cal{M}}M\neq
M$. Thus by Proposition 2.9, $M$ is a ${\mathbb{M}}$-radical
module.  $\square$

\begin{Pro}  Let $M$ be a semisimple $R$-module.  Then the following statement are equivalent:\vspace{3mm}\\
{\rm  (1)} $M$ is a ${\mathbb{P}}$-radical module.\vspace{1mm}\\
 {\rm  (2)} $M$ is a ${\mathbb{M}}$-radical module and $R/Ann(M)$ is a Hilbert ring.\vspace{1mm}\\
  {\rm  (3)} $M$ is full semisimple and $R/Ann(M)$ is a Hilbert ring.
\end{Pro}
\noindent {\bf Proof.} Since $M$ is a semisimple $R$-module, we
can assume that $M=\bigoplus_{i\in I}R/{\cal{M}}_i$ where $I$ is
an index set and each  ${\cal{M}}_i$ is a maximal ideal of $R$.\\
 $(1)\Rightarrow(2)$. Since every ${\mathbb{P}}$-radical module is ${\mathbb{M}}$-radical,
  it is suffices to show that $R/{\rm Ann}(M)$ is a Hilbert ring. Suppose ${\cal{P}}\supseteq {\rm Ann}(M)$ is a prime
  ideal of $R$. Thus $(\sqrt[p]{{\cal{P}}M}:M)={\cal{P}}$ and hence
    $\sqrt[p]{\cal{P}M}\neq M$. We can assume that
$\sqrt[p]{\cal{P}M}=\bigcap_{\lambda\in\Lambda}P_\lambda$, where
$\Lambda$ is an index set and each $P_\lambda$ is a prime
submodule of $M$ containing  ${\cal{P}M}$. Thus for each
$P_\lambda$, the factor module  $M/P_\lambda$ is a prime
     semisimple module and hence by Lemma 3.4, $M/P_\lambda$ is a homogenous semisimple $R$-module, i.e., ${\cal{M}}_\lambda:=(P_\lambda:M)$ is a
     maximal ideal of $R$. Thus $${\cal{P}}=(\sqrt[p]{{\cal{P}}M}:M)=(\bigcap_{\lambda\in\Lambda}P_\lambda:M)=
     \bigcap_{\lambda\in\Lambda}(P_\lambda:M)=\bigcap_{\lambda\in\Lambda}{\cal{M}}_\lambda.$$
     Thus every prime ideal ${\cal{P}}\supseteq {\rm Ann}(M)$ is an
     intersection of maximal ideals of $R$, i.e., $R/{\rm Ann}(M)$ is a Hilbert ring.\\
      $(2)\Leftrightarrow(3)$ is by Proposition 3.5.\\
      $(2)\Rightarrow(1)$ is by Theorem 2.11.  $\square$
\begin{Cor}  Let $R$ be a Hilbert ring or a domain with $dim(R)=1$. For a  semisimple $R$-module $M$ the following statement are equivalent:\vspace{3mm}\\
{\rm  (1)} $M$ is a ${\mathbb{P}}$-radical module.\vspace{1mm}\\
 {\rm  (2)} $M$ is a ${\mathbb{M}}$-radical module.\vspace{1mm}\\
  {\rm  (3)} $M$ is full semisimple.
\end{Cor}
\noindent {\bf Proof.} If $R$ is a Hilbert ring then by
Proposition 3.5 and Theorem 2.11,  the proof is complete. Thus we
can assume that $R$ is a domain with $dim(R)=1$. Since $M$ is a
semisimple $R$-module, we can assume that $M=\bigoplus_{i\in
I}R/{\cal{M}}_i$ where $I$ is
an index set and each  ${\cal{M}}_i$ is a maximal ideal of $R$.\\
 $(1)\Rightarrow(2)$ is clear.\\
$(2)\Rightarrow(3)$ is by Proposition 3.5.\\
 $(3)\Rightarrow(1)$.
Suppose ${\cal{P}}\supseteq {\rm Ann}(M)$ is a prime  ideal of
$R$. Let ${\cal{M}}\supseteq {\cal{P}}$ be a maximal ideal. Since
$R$ is a domain with $dim(R)=1$, either ${\cal{P}}=(0)$ or
${\cal{M}}={\cal{P}}$. If ${\cal{M}}={\cal{P}}$, then
${\cal{P}}={\cal{M}}_i$ for some $i$ and so
${\cal{P}}M={\cal{M}}_iM\neq M$. Clearly  ${\cal{P}}M$ is a prime
submodule of $M$ with $({\cal{P}}M:M)={\cal{P}}$. It follows that
$(\sqrt[p]{{\cal{P}}M}:M)=({\cal{P}}M:M)={\cal{P}}$. Now assume
that  ${\cal{P}}=(0)$. Then ${\rm Ann}(M)=(0)=\bigcap_{i\in
I}{\cal{M}}_i$. Since every proper submodule of a simisimple
module is an intersection of maximal submodules and each maximal
submodule is a prime submodule, we conclude that
$\sqrt[p]{\cal{P}M}=\sqrt[p]{(0)}=(0)$. It follows that
$(\sqrt[p]{{(0)}M}:M)=((0):M)={\rm Ann}(M)=(0)=\sqrt{(0)}$. Thus  $M$ is a ${\mathbb{P}}$-radical module.  $\square$\\

We conclude this paper   with the following result that offer
several characterizations for semisimple primeful modules.

\begin{Cor}  Let $M$ be a semisimple $R$-module.  Then the following statement are equivalent:\vspace{3mm}\\
{\rm  (1)} $M$ is a primeful  module.\vspace{1mm}\\
{\rm  (2)} $M$ is a ${\mathbb{P}}$-radical module and $dim(R/Ann(M))=0$.\vspace{1mm}\\
 {\rm  (3)} $M$ is a ${\mathbb{M}}$-radical module and $dim(R/Ann(M))=0$.\vspace{1mm}\\
  {\rm  (4)} $M$ is a full semisimple  module and $dim(R/Ann(M))=0$.
\end{Cor}
\noindent {\bf Proof.} $(1)\Rightarrow(2)$.  Since every primful
module is ${\mathbb{P}}$-radical,
  it is suffices to show that $dim(R/{\rm Ann}(M))=0$. Suppose
${\cal{P}}\supseteq {\rm Ann}(M)$ is a prime  ideal of $R$. Thus
there exists  a prime submodule $P$ of $M$ such that
$(P:M)={\cal{P}}$. Since $M/P$ is a prime semisimple $R$-module,
by Lemma 3.4,
${\cal{P}}$ is a maximal ideal. Thus $dim(R/{\rm Ann}(M))=0$.\\
$(2)\Rightarrow(3)$ is by by Theorem 2.12.\\
$(3)\Rightarrow(4)$ is by by Proposition 3.5. $\square$\\

\noindent {\bf Note:} In Part II we shall continue the study of
this construction.

\end{document}